\documentclass{icmart}

\newtheorem{theorem}{Theorem}[section]
\newtheorem{lemma}[theorem]{Lemma}

\theoremstyle{definition}
\newtheorem{definition}[theorem]{Definition}

\newtheorem{conjecture}[theorem]{Conjecture}

\newcommand{\ZZ}{\mathbb{Z}}

\newcommand{\QQ}{\mathbb{Q}}
\newcommand{\CC}{\mathbb{C}}
\newcommand{\RR}{\mathbb{R}}

\def\T{\Theta}
\def\a{\alpha}

\def\D{\Delta}

\def\g{\gamma}

\def\G{\Gamma}

\def\l{\lambda}

\def\eg{{\it e.g.}\ }
\def\ie{{\it i.e.}\ }

\def\PSL{\mbox{\rm{PSL}}}

\def\Isom{\mbox{Isom}}
\def\Mob{\mbox{M\"ob}}
\def\PGL{\mbox{\rm{PGL}}}

\def\GL{\mbox{\rm{GL}}}

\def\qed{ $\Box$}

\def\OO{\mathcal{O}}

\def\Vol{\mbox{\rm{Vol}}}
\def\Comm{\mbox{\rm{Comm}}}

\def\HH{\mathbb{H}}
\def\SS{\mathbb{S}}

\def\PP{\mathbb{P}}
\def\CC{\mathbb{C}}
\def\EE{\mathbb{E}}

\def\interior{\rm int}

\edef\t@mp{\catcode`\noexpand\~=\the\catcode`\~}%
    \catcode`\~=12
    \def\tild@{~}%
    \t@mp

\title{Finiteness of arithmetic Kleinian reflection groups} 

\author[Ian Agol]{%
        Ian Agol \thanks{Agol partially supported by NSF grant DMS-0204142 and the Sloan Foundation}
} 

\contact[agol@math.uic.edu]{MSCS UIC 322 SEO, m/c 249\\
        851 S. Morgan St.\\
        Chicago, IL 60607-7045}

\date{%
 December 21, 2005}

\begin{document}

\begin{abstract} 
We prove that there are only finitely many arithmetic Kleinian maximal reflection
groups. 
\end{abstract} 



\maketitle

\section{Introduction}
A Kleinian reflection group $\G$ is a discrete group generated by reflections in the 
faces of hyperbolic polyhedron $P\subset \HH^{3}$. We may assume that the 
dihedral angles of $P$ are of the form $\pi/n$, $n\geq 2$, in which case $P$ forms
a fundamental domain for the action of $\G$ on $\HH^{3}$. If $P$ has finite volume,
then $\HH^{3}/\G=\mathcal{O}$ is a hyperbolic orbifold of finite volume, which is
obtained by ``mirroring'' the faces of $P$. Andreev gave a combinatorial
characterization of hyperbolic reflection groups in 3-dimensions, in terms of the
topological type of $P$ and the dihedral angles assigned to the edges of $P$ \cite{Andreev70}. 
A reflection group $\G$ is maximal if there is no reflection group $\G'$ such that
$\G < \G'$. We'll defer the definition of arithmetic
groups until later, but a theorem of Margulis implies that $\G$ is arithmetic if and only if $[\Comm(\G):\G]=\infty$,
where $\Comm(\G)=\{g\in \Isom(\HH^{3}) | [\G:g^{-1}\G g \cap \G]<\infty\}$ \cite{Margulis75}.
The main theorem of this paper is that there are only finitely many arithmetic Kleinian
groups which are maximal reflection groups.

The argument generalizes an argument of Long-MacLachlan-Reid \cite{LMR05}, which implies that
there are only finitely many arithmetic minimal hyperbolic 2-orbifolds with bounded
genus. Their argument is in fact a generalization of an argument of Zograf \cite{Zograf91}, who 
reproved that there are only finitely many congruence groups $\G$ commensurable with
$\PSL(2,\ZZ)$ such that $\HH^{2}/\G$ has genus $0$ 
(this was proven originally by Dennin \cite{Dennin71,Dennin72}, and was known as
Rademacher's conjecture). 
The key ingredient of their argument is a theorem of Vigneras \cite{Vigneras83} (based on work
of Gelbart-Jacquet  \cite{GJ78} and Jacquet-Langlands \cite{JL70}), which implies that a congruence arithmetic
Fuchsian group has $\l_{1}\geq \frac{3}{16}$ (and which has a generalization to higher
dimensions \cite{BurgerSarnak91}). The other key ingredient is an estimate of Zograf \cite{Zograf84}, which
implies that for a hyperbolic 2-orbifold $\mathcal{O}$, $\l_{1}(\mathcal{O})\Vol(\mathcal{O})$ is bounded
linearly by the genus of $\mathcal{O}$. Zograf's argument is based on a sequence of
improvements of a result of Sz\"{e}go (who did this for planar domains)\cite{Szego54}, 
by Hersch (for $S^{2}$)\cite{Hersch70}, 
Yang-Yau (for orientable surfaces) \cite{YangYau80}, and Li-Yau (for non-orientable surfaces and manifolds
of a fixed conformal type)\cite{LiYau82}. We observe that the Li-Yau argument (sharpened by El Soufi-Ilias \cite{ElSoufiIlias86})
generalizes to orbifolds, and we then apply this to arithmetic Kleinian reflection groups. 

In the concluding section, we consider how one might generalize this result to higher
dimensions to prove 
\begin{conjecture} \label{finite}
There are only finitely many maximal arithmetic reflection groups in $\Isom(\HH^{n})$,
$n>1$.
\end{conjecture}

\section{Conformal volume of orbifolds}
Conformal volume was first defined by Li-Yau, partially motivated by generalizing results
on surfaces due to Yang-Yau, Hersch, and Sz\"ego. We generalize this notion to orbifolds. 
Let $(\mathcal{O},g)$ be a compact Riemannian orbifold, possibly with boundary. Let
$|\mathcal{O}|$ denote the underlying topological space. Denote 
the volume form by $dv_{g}$, and $\Vol(\mathcal{O},g)$ its volume. Let $\Mob(\SS^{n})$ denote
the conformal transformations of $\SS^{n}$. It is well-known that $\Mob(\SS^{n})=\Isom(\HH^{n+1})$. 
$|\mathcal{O}|$ has a codimension zero dense open subset which is a Riemannian manifold. We will
say that a map 
$\varphi: |\mathcal{O}_{1}|  \to |\mathcal{O}_{2}|$ is {\it PC} if it is a continous map which is piecewise conformal
on the open submanifold of $|\mathcal{O}_{1}|$  which maps to the manifold part of $|\mathcal{O}_{2}|$. 
Clearly, if $\varphi:|\mathcal{O}|\to \SS^{n}$ is PC, and $\mu\in \Mob(\SS^{n})$, then $\mu\circ \varphi$ is also PC. 

\begin{definition}
For a piecewise smooth map $\varphi: |\mathcal{O}| \to (\SS^{n},can)$, define

$$V_{c}(n,\varphi) = \underset{\mu\in \Mob(\SS^{n})}{\sup} \Vol(\mathcal{O}, (\mu\circ \varphi)^{*}(can)).$$
If there exists a PC map $\varphi: |\mathcal{O}| \to \SS^{n}$, then we also define
$$V_{c}(n, \mathcal{O}) = \underset{\varphi:\mathcal{O}\to \SS^{n}  PC}{\inf} V_{c}(n,\varphi).$$

$V_{c}(n,\mathcal{O})$ is denoted the ($n$-dimensional) {\it conformal volume} of $\mathcal{O}$. 
\end{definition}

{\bf Remark:} For our application, it would suffice  to define a PC map to be Lipschitz and {\it a.e.} conformal. 
It seems likely that our definition of conformal volume coincides with that of Li-Yau for manifolds,
but we have not checked this (it would suffice to show that a PC map can be approximated
by conformal maps). 

If there exists a piecewise isometric map $\varphi: (|\mathcal{O}|,g) \to \EE^{n}$ for some $n$,
then clearly $V_{c}(n,\mathcal{O})$ is well-defined, since $\EE^{n}$ has a conformal embedding
into $\SS^{n}$. For an orbifold $(\mathcal{O},g)$, to prove that $V_{c}(n,\mathcal{O})$ is well-defined, 
one need only check that  $(|\mathcal{O}|,g)$ embeds piecewise isometrically into
some compact Riemannian manifold, in which case the Nash-Moser isometric embedding theorem 
implies that $(|\mathcal{O}|,g)$ embeds piecewise isometrically into $\EE^{n}$, for some $n$. 
We will only apply conformal volume to orbifolds which will obviously have a  PC map to $\SS^{n}$,
for some $n$. We record basic facts about conformal volume which were recorded in \cite{LiYau82},
and which carry over to our notion of conformal volume for orbifolds.

Fact 1: If $\mathcal{O}$ admits a degree $d$ PC map onto another orbifold $\mathcal{P}$, then
$$V_{c}(n,\mathcal{O}) \leq |d| V_{c}(n,\mathcal{P}).$$ 

Fact 2: Since $\SS^{n}\hookrightarrow \SS^{n+1}$ embeds isometrically, it's clear that $V_{c}(n,\mathcal{O})\geq V_{c}(n+1,\mathcal{O})$. 
Define the conformal volume $V_{c}(\mathcal{O})=\underset{n\to \infty}{\lim} V_{c}(n,\mathcal{O})$. 

Fact 3: If $\mathcal{O}$ is of dimension $m$, and $\varphi: |\mathcal{O}| \to \SS^{n}$ is a PC map, then 
$$V_{c}(n,\varphi)\geq V_{c}(n,\SS^{m})=\Vol(\SS^{m}).$$ 
The same argument as in Li-Yau works here: ``blow up'' about a smooth  manifold point of $\varphi(|\mathcal{O}|)$
so that  the image Hausdorff limits to a geodesic sphere of dimension $m$.

Fact 4: If $\mathcal{O}$ is an embedded suborbifold of the orbifold $\mathcal{P}$, and $\varphi: |\mathcal{P}|\to \SS^{n}$ 
is  PC, then $V_{c}(n, \varphi) \geq V_{c}(n,\varphi|_{|\mathcal{O}|})$.
Thus, $V_{c}(n,\mathcal{P}) \geq V_{c}(n, \mathcal{O})$.

\section{Finite subgroups of $O(3)$}

\begin{lemma}
Let $G<O(3)$ be a finite subgroup. Then there exists a group $G'$, $G\leq G'<O(3)$, which is 
generated by reflections, such that $[G':G]\leq 4$. 
\end{lemma}
\begin{proof}
This follows from a case-by case analysis of spherical 2-orbifolds. 
We will use Conway's notation for spherical
orbifolds (see \eg  ch. 13 \cite{Th}). 

Case $(*), (*p,p)$ or $(*p,q,r)$: These orbifold groups are generated by reflections.

Case $(p,p)$ or $(p,q,r)$: These 2-fold cover $(*p,p)$ or $(*p,q,r)$ respectively.

Case $(n *)$:  This 2-fold covers  $(*n,2,2)$.

Case $(2*m)$: This 2-fold covers $(*2m,2,2)$.

Case $(3*2)$: This 2-fold covers $(*4,3,2)$.

Case $(n | \circ)$: This 2-fold covers $(2n *) \overset{2:1}{\to} (* 2n,2,2) $ (this includes
the case $n=1$, \ie  $(1|\circ)= \RR\PP^{2}$). 

This exhausts all possible spherical orbifolds, and we see that in each case the orbifold
fundamental group  is  of index $\leq 4$ in a reflection group (all but the last case have
index $\leq 2$).
\end{proof}

{\it Question:} Given a dimension $n$, is there a constant $C(n)$ such that any finite
subgroup of $O(n)$ is of index $\leq C(n)$ in a reflection group? If so, then one should 
be able to prove conjecture \ref{ellipticconformalvolume}. We suspect
the answer to this question is no, which is why we have not been
able to generalize the proof in this section to higher dimensions. 

\section{Eigenvalue bounds}

We observe that the argument of Thm. 2.2 \cite{ElSoufiIlias86} generalizes to our
context of orbifolds (their theorem sharpens Cor. 3 sect. 2 of \cite{LiYau82}). If $(\mathcal{O},g)$
is a Riemannian orbifold with piecewise smooth boundary, then $\l_{1}(\mathcal{O},g)$ 
is the first non-zero eigenvalue of $\D_{g}$ on $\mathcal{O}$ with Neumann boundary
conditions. 

\begin{theorem} \cite{ElSoufiIlias86} \label{ev}
Let $(\mathcal{O},g)$ be a compact Riemannian orbifold of dimension $m$, possibly
with boundary. If $\varphi:|\mathcal{O}|\to \SS^{n}$ is a PC map, 
$$\l_{1}(\mathcal{O},g) \Vol(\mathcal{O},g)^{\frac2m} \leq m V_{c}(n,\varphi)^{\frac2m}.$$
\end{theorem}
\begin{proof}
Let ${\bf X} = (X_{1}, ..., X_{n+1})$ be the coordinate
functions on $\SS^{n}\subset \RR^{n+1}$. Then $\sum_{i=1}^{n+1} X_{i}^{2}= 1$, restricted to $\SS^{n}$.

\begin{lemma}
There exists $\mu \in \Mob(\SS^{n})$ such that 
$\int_{\mathcal{O}} {\bf X}\circ \mu \circ \varphi\ dv_{g}={\bf 0}$.
\end{lemma}

\begin{proof}
Let ${\bf x} \in \SS^{n}$. For $0\leq t<1$, $t {\bf x}\in \HH^{n+1}$,
let $\mu_{t{\bf x}} \in \Mob(\SS^{n}) = \Isom(\HH^{n+1})$
be the hyperbolic translation along the ray $\RR{\bf x}$ taking ${\bf 0}$
to $t{\bf x}$ (thus, $\mu_{0 {\bf  x}}=\mu_{{\bf 0}}=Id$). 
Let $H(t{\bf x}) = \frac{1}{\Vol(\mathcal{O},g)} \int_{\mathcal{O}} {\bf X}\circ \mu_{t{\bf x}}\circ \varphi\ dv_{g}$.
We may think of $H$ as defining a function $H:\HH^{n+1}\to \HH^{n+1}$, 
which gives the center of mass of the measure
coming from $\varphi_{*} dv_{g}$, where
we take the point $-t{\bf x}$ to be the origin of the sphere
$\SS^{n}=\partial \HH^{n+1}$ by the conformal map $\mu_{t{\bf x}}$. As $t\to 1$,
all of the mass of $\mu_{t\bf{x}} \circ \varphi (\mathcal{O})$ becomes
concentrated at ${\bf x}$, and we see that $H$ extends to a
continuous function $H:B^{n+1}\to B^{n+1}$ (where $B^{n+1}=\HH^{n+1}\cup \SS^{n}$) such that
$H_{|\SS^{n}}=Id_{|\SS^{n}}$.  Thus, $H$ is onto, so there exists ${\bf y} \in \HH^{n+1}$
such that $H({\bf y}) = {\bf 0}$, and we take $\mu=\mu_{{\bf y}}$. \end{proof}

Now, replace $\varphi$ with $\mu\circ \varphi$, noting that this
is still a PC map. Then $X_{i}\circ \varphi$ may be 
used as test functions in the Rayleigh-Ritz quotient, since they
are Lipschitz functions which are $L^{2}$ orthogonal to the
constant function. Thus,
$$\l_{1}(\mathcal{O}) \int_{\mathcal{O}} |X_{i}\circ \varphi|^{2} dv_{g}\leq \int_{\mathcal{O}} |\nabla X_{i}\circ \varphi |^{2} dv_{g}.$$
Summing, we see that 
$$\l_{1}(\mathcal{O}) \int_{\mathcal{O}} \sum_{i=1}^{n+1}   |X_{i}\circ \varphi|^{2} dv_{g} = \l_{1}(\mathcal{O}) \Vol(\mathcal{O},g) 
\leq m \int_{\mathcal{O}} \frac1m \sum_{i=1}^{n+1} |\nabla X_{i}\circ \varphi |^{2} dv_{g}$$ 
$$ \leq m \left(\int_{\mathcal{O}} \left( \frac1m \sum_{i=1}^{n+1} |\nabla X_{i}\circ \varphi|^{2}\right)^{\frac{m}{2}} dv_{g}\right)^{\frac2m} \Vol(\mathcal{O},g)^{1-\frac2m} ,$$
where the last inequality is H\"olders inequality. 
Now, we use the fact that $\varphi$ is PC to see that $\varphi^{*} can = \frac1m \sum_{i=1}^{n+1} |\nabla X_{i}\circ \varphi|^{2} a.e.$,
and thus 
$$\int_{\mathcal{O}} \left( \frac1m \sum_{i=1}^{n+1} |\nabla X_{i} \circ \varphi |^{2} \right) ^{\frac{m}{2}} dv_{g}  \leq \Vol(\mathcal{O},\varphi^{*}can).$$
Finally, we obtain the desired inequality putting these inequalities together. 
\end{proof}

\section{Congruence arithmetic hyperbolic 3-orbifolds}
We need to know some properties of arithmetic hyperbolic 3-orbifolds. For background
and notation, see Maclachlan-Reid \cite{MR03}.

Let $k\subset \CC$ be a number field with only one complex place, and let $A$ be a quaternion
algebra over $k$. Let $\rho: A\to M(2,\CC)$ be a $k$-embedding, $P:\GL(2,\CC)\to \PGL(2,\CC)$. Let $\mathcal{E}\subset A$ be
either a maximal order or an Eichler order, and let $N(\mathcal{E})\subset A^{*}$ be the normalizer of 
$\mathcal{E}$ in $A$. Let $\G$ be an arithmetic Kleinian group, such
that $\G=P\rho G$, $G \subset A^{*}$. Then there exists an
order $\mathcal{E}\subset A$ which is either a maximal order or an Eichler order such that 
$G\subset N(\mathcal{E})$ (see thm. 11.4.3 \cite{MR03}). In particular, if $\G$ 
is a maximal arithmetic Kleinian group, then $\G=P\rho N(\mathcal{E})$, for some order $\mathcal{E}$. 

An ideal $I$ in $A$ is a complete $R_{k}$-lattice. 
The left order of $I$ is $\mathcal{O}_{l}(I)= \{ a \in A | aI \subset I \}$, and the
right order is $\mathcal{O}_{r}(I) = \{ a \in A | Ia\subset I\}$. The ideal $I $ is {\it 2-sided}
if $\mathcal{O}_{l}(I)=\mathcal{O}_{r}(I)$. The ideal is {\it integral} if 
$I$ lies in both $\mathcal{O}_{l}(I)$ and in $\mathcal{O}_{r}(I)$ (\ie $I^{2}\subseteq I$ is multiplicatively closed). 
If $\mathcal{O}$ is a maximal order in $A$, and
$I$ is a 2-sided integral ideal in $\mathcal{O}$ such that $\mathcal{O}=\mathcal{O}_{l}(I)=\mathcal{O}_{r}(I)$, 
then the \textit{principal congruence subgroup} of $\mathcal{O}^{1}$ is 
$$\mathcal{O}^{1}(I) = \mathcal{O}^{1} \cap (1+I).$$

Thus, $\mathcal{O}^{1}(I)$ is the kernel of the map $\mathcal{O}^{1} \to \mathcal{O}/I$, which is
therefore of finite index in $\mathcal{O}^{1}$, since $\mathcal{O}/I$ is finite. A discrete group $G \subset A$ is {\it congruence} if 
it contains a principal congruence subgroup, and $\G<\PGL(2,\CC)$ is congruence if $\G=P\rho G$, for
some $G<A$ congruence. 

\begin{lemma} (Long-MacLachlan-Reid \cite{LMR05})
A maximal arithmetic Kleinian group is congruence. 
\end{lemma}
\begin{proof}
Let $\G\subset PGL(2,\CC)$ be a maximal arithmetic Kleinian group. Then $\G=P\rho N(\mathcal{E})$, for
some order $\mathcal{E}$ of square-free level.  If $\mathcal{E}$ is a maximal order, then $\mathcal{E}^{1}$
is a congruence subgroup for the trivial ideal $I=\mathcal{E}$. If $\mathcal{E} =\mathcal{O}_{1}\cap\mathcal{O}_{2}$,
where $\mathcal{O}_{i}$ are maximal orders (so that $\mathcal{E}$ is an Eichler order), then choose
$\a \in R_{k}-\{0\}$ such that $I=\a \mathcal{O}_{1} \subset \mathcal{O}_{2}$. Then $\mathcal{O}_{l}(I)=\mathcal{O}_{r}(I)=\mathcal{O}_{1}$, and
$I^{2}=\a^{2} \mathcal{O}_{1} \subset \a \mathcal{O}_{1}=I$. Also, clearly $1+I \subset \OO_{1}\cap \OO_{2}$. Thus, $\mathcal{O}^{1}(I) =\mathcal{O}_{1}^{1}\cap (1+I) \subset \OO_{1}\cap \OO_{2}$.
Thus, $\mathcal{O}^{1}(I)\subset \mathcal{E}^{1}$,   and we see that $\mathcal{E}^{1}$ is a congruence subgroup, and therefore
$P\rho N(\mathcal{E})$ is congruence. 
\end{proof}

A fundamental theorem of Vigneras, making use of work of Jacquet-Langlands \cite{JL70} and Gelbart-Jacquet \cite{GJ78},
generalizes a result of Selberg for $\PSL(2,\ZZ)$. For a hyperbolic orbifold $\OO$, let $\l_{1}(\OO)$
be the minimal non-zero eigenvalue of the Laplacian $\D$ on $\OO$.

\begin{theorem} \cite{Vigneras83, BurgerSarnak91}
Let $\OO=\HH^{3}/\G$, where $\G$ is an arithmetic congruence subgroup. Then $\l_{1}(\OO) \geq \frac34$. 
\end{theorem} 

It is conjectured that under the hypotheses of the above theorem, $\l_{1}(\OO)\geq 1$, which
is known as (a special case of) the {\it generalized Ramanujan conjecture} \cite{BurgerSarnak91}.

\section{Finiteness of arithmetic Kleinian maximal reflection groups}

We put together the results in the previous sections to prove our main theorem.

\begin{theorem}
There are only finitely many arithmetic maximal reflection groups in $\Isom(\HH^{3})$. 
\end{theorem}
\begin{proof}
Suppose that $\G$ is an arithmetic maximal reflection group. That is, there
is no group $\G' < \Isom(\HH^{3})$, with $\G <\G'$, such that $\G'$ is generated by reflections.
Then there exists $\G\leq\G_{0}<\Isom(\HH^{3})$, such that $\G_{0}$ is a maximal Kleinian group.
$\G$ is generated by reflections in a finite volume polyhedron $P$. 

\begin{lemma} (Vinberg \cite{Vinberg67})
$\G$ is a normal subgroup of $\G_{0}$. Moreover, there is a finite subgroup $\T<\G_{0}$
such that $\T \to  \G_{0}/\G$ is an isomorphism, and $\T$ preserves the polyhedron $P$. 
\end{lemma}
\begin{proof}
This  follows from the fact that the set of reflections in $\G_{0}$
is conjugacy invariant, and therefore the group generated by reflections is normal
in $\G_{0}$. Since $\G$ is a maximal reflection group, this subgroup must be $\G$.
The polyhedron $P$ is a fundamental domain of $\G$, and if there is an element
$\g\in \G_{0}$ such that $\interior(P)\cap \g(\interior(P)) \neq \emptyset$, then $\g(P)=P$.
Otherwise, there would be a geodesic plane $V$ containing a face of $P$, such
that $V\cap \interior(P) \neq \emptyset$. The reflection $r_{V}\in \G$ in the plane $V$
would be conjugate to a reflection $r_{\g(V)}= \g r_{V} \g^{-1}$, which is not
in $\G$ since $r_{g(V)} (\interior(P))\cap \interior(P) \neq \emptyset$, which implies that
$\G$ is not a maximal reflection group, a contradiction. Let $\T$ be the subgroup of
$\G_{0}$ such that $\T(P)=P$. Clearly $\T$ is finite, since $P$ is finite volume and
has finitely many faces. If
$\g_{0}\in \G_{0}$, let $\g\in \G$ be such that $\g_{0}(\interior(P)) \cap \g(\interior(P))\neq \emptyset$.
Then $\g^{-1} \g_{0} (P)=P$, so $\g_{0} \in \g \T$. Thus, $\T\to \G_{0}/\G$ is
an isomorphism. 
\end{proof}
Let $\mathcal{O}=\HH^{3}/\G_{0}$, and $\T$ is the finite group coming from the previous lemma.

\begin{lemma}
$\l_{1}(\mathcal{O}) = \l_{1}(P/\T)$.
\end{lemma}
\begin{proof}
Let $f$ be an eigenfunction on $P/\T$ with eigenvalue $\l_{1}(P/\T)$. Since
$f$ has Neumann boundary conditions, its level sets in $P/\T$ are orthogonal
to $\partial P/\T$. Let $\tilde{f}$ be the preimage of $f$ under the map $P\to P/\T$,
so that $\tilde{f}$ is invariant under the action of $\T$. Extend $\tilde{f}$ to
a function $\tilde{F}$ on $\HH^{3}$, by the action of $\G$ (and therefore invariant
 under the action of $\G_{0}$). By the reflection principle, $\tilde{F}$ is a
smooth function, so it descends to an eigenfunction $F$ of $\D$ on $\mathcal{O}$
with eigenvalue $\l_{1}(P/\T)$. Conversely, if $F$ is an eigenfunction of 
$\D$ on $\mathcal{O}$ with eigenvalue $\l_{1}(\mathcal{O})$, then 
$F_{| P/\T}$ gives an eigenfunction on $P/\T$ with Neumann boundary conditions,
since the level sets of $\tilde{F}$ must be invariant under reflections, and therefore
perpendicular to the faces of $P$. 
\end{proof}

Consider $\HH^{3}\subset \SS^{3}$ embedded conformally as the upper half
space of $\SS^{3}$, so that $\Isom(\HH^{3})$ acts conformally on $\SS^{3}$. 
Normalize so that $\T$ acts isometrically on $\SS^{3}$, which we may
do by a result of Wilker \cite{Wilker82}. Clearly, $V_{c}(n,\mathcal{O})=V_{c}(n,P/\T)$,
since $|\OO|=|P/\T|$. 
Then the orbifold $P/\T \subset \SS^{3}/\T$ is a conformal embedding, so by
fact 4, $V_{c}(3,P/\T) \leq V_{c}(3,\SS^{3}/\T)$. The group $\T$ embeds in a finite
reflection group $\T' \subset O(3)$ such that $[\T':\T]\leq 4$. Then by fact  1, 
$V_{c}(3,\SS^{3}/\T) \leq 4 V_{c}(3,\SS^{3}/\T')$. Now, there is a polyhedron
$Q\subset \SS^{3}$ with geodesic faces which is the fundamental domain for
$\T'$. Clearly $V_{c}(3,Q)= V_{c}(3,\SS^{3}/\T')$. By facts 2 and 4, $V_{c}(3,Q)=\Vol(S^{3}, can)=2\pi^{2}$.
Thus, we have $V_{c}(\mathcal{O}) \leq 8\pi^{2}$. 

Now, we apply the eigenvalue estimates:
$$\frac34 \Vol(\mathcal{O}  )^{\frac23} \leq \l_{1}(\mathcal{O}) \Vol(\mathcal{O})^{\frac23}  \leq 3 V_{c}(\HH^{3}/\G_{0})^{\frac23}\leq 3 (8\pi^{2})^{\frac23}.$$
Thus, we obtain $\Vol(\mathcal{O}) \leq 64\pi^{2}.$ Since volumes of arithmetic hyperbolic orbifolds
are  discrete, and $\G_{0}$ is uniquely determined by $\G$, we conclude that there are only finitely
many arithmetic maximal reflection groups. 
\end{proof}

\section{Conclusion}

From the main theorem, we conclude that given an arithmetic reflection group  $\G<\Isom(\HH^{3})$, it must lie in 
one of finitely many maximal reflection groups (up to conjugacy). If $\G$ is a reflection group in a polyhedron $P$ for which all the
dihedral angles are $\pi/2$ or all are $\pi/3$, then
there are infinitely many co-finite volume reflection subgroups of $\G$. Thus we see
that there are commensurability classes of arithmetic groups for which there are infinitely many
reflection groups in the commensurability class, and thus in our finiteness result, the
maximality assumption is crucial. It is an interesting project
to try to identify all arithmetic maximal reflection groups, and to classify their reflection subgroups
of finite covolume.

 For the non-compact examples, one
may apply volume formulae to estimate the maximal
discriminant of a quadratic imaginary number field $k$ for
which $\PGL(2,k)$ contains a reflection group. 
Humbert first computed the covolumes of Bianchi groups, and a generalization due to 
Borel implies that for a non-compact arithmetic Kleinian group, the minimal covolume $\mu$
satisfies 
$$\mu \geq \frac{|\D_{k}|^{\frac32} \zeta_{k}(2)}{16\pi^{2} h_{k}},$$ where $h_{k}$ is the 
class number of the number field $k$. The Brauer-Siegel theorem gives an estimate of 
$h_{k}$ for a number field, and implies for a quadratic imaginary 
number field $k$ that 
$$ |\D_{k}| \zeta_{k}(2) \geq \frac{h_{k}(2\pi)^{2}}{2w},$$
where $w$ is the order of the group of roots of unity in $k$. If $k\neq \QQ(i), \QQ(\sqrt{-3})$,
then $w=2$ (see the proof of theorem 11.7.2 \cite{MR03}). Thus, we have 
$$64\pi^{2}\geq \mu \geq \frac{1}{16} |\D_{k}|^{\frac12}.$$  
Thus 
$|\D_{k}| < 2^{20} \pi^{4}=1.02\times 10^{8}$. Hatcher has computed orbifold structures of 
some of the Bianchi groups of small discriminant, and from his pictures one may deduce that
$\PGL(2,R_{k})$ is commensurable with a reflection group for 
$$\D_{k}= -3,-4,-7,-8,-11,-15,-19,-20,-24, -39, -40,-52,-55,-56,-68,-84,$$ where $k$ is a quadratic imaginary number field \cite{Hatcher83}. In principle, it ought to 
be possible to compute all arithmetic reflection groups, but clearly even in the non-compact
case, the volume estimates we obtain do not make this computation feasible. To classify
 non-compact Kleinian arithmetic groups, it may require the
infusion of some more number theory.  Arithmetic restrictions have been found on reflection
groups containing the Bianchi groups by Vinberg \cite{Vinberg90} and Blume-Nienhaus \cite{B-N92}.
If these results could be extended to all maximal non-compact
arithmetic groups, then one may be able to give a complete classification. The
classification of compact arithmetic maximal reflection groups, although in principle decidable,
is probably not feasible at this stage.

It's clear that for a finite volume polyhedron 
$P\subset \HH^{n}$, $V_{c}(n,P)= \Vol(\SS^{n})$, so 
$$\l_{1}(P) \Vol(P)^{\frac2n} < n \Vol(\SS^{n})^{\frac2n}.$$
Thus  in $n$ dimensions, there
are finitely many reflection groups $\G< \Isom(\HH^{n})$ which have a lower bound on 
$\l_{1}(\HH^{n}/\G)$. It would be interesting
if one could generalize the arguments of the main theorem to higher dimensions. It is known by
work of Prokhorov that there cannot be any reflection groups in dimension $>996$ \cite{Prokhorov86}
(Vinberg showed that compact reflection groups cannot exist in dimension $>30$ \cite{Vinberg84}). 
Vinberg also gave a characterization of arithmetic reflection groups, in terms of a totally 
real field of definition $K$ and an integrality condition \cite{Vinberg67}. 
Nikulin has shown that there exists a constant $N_{0}$ such that the
set of arithmetic groups in $\Isom(\HH^{n})$ generated by reflections with   $n>16$ and $[K:\QQ] > N_{0}$ is
empty \cite{Nikulin81}. It was proved in a previous paper by  Nikulin that the set of maximal arithmetic groups
generated by reflections in $\Isom(\HH^{n})$ with a fixed degree $[K:\QQ]$
is finite \cite{Nikulin80}.
Thus,
to generalize the main argument of this paper, 
one would have to bound the conformal volume of a minimal arithmetic $n$-orbifold containing a
reflection suborbifold up to dimensions $n\leq 16$.  By the characterization of maximal arithmetic
groups in $\Isom(\HH^{n})$, it should follow that they are congruence. By a result of Burger-Sarnak,
if $\G< \Isom(\HH^{n})$ is a congruence arithmetic reflection group with $n>2$, then 
$\l_{1}(\HH^{n}/\G) > \frac{2n-3}{4}$ (Cor. 1.3 \cite{BurgerSarnak91}, and the fact
that reflection groups are defined by quadratic forms). Every maximal arithmetic
group covered by a reflection group will embed conformally in an elliptic $n$-orbifold. So to
to prove the following conjecture for $n\leq 16$:

\begin{conjecture} \label{ellipticconformalvolume}
There is a function $K(n)$, such that if $\mathcal{O}$ is an elliptic $n$-orbifold, 
then $V_{c}(\mathcal{O})\leq K(n)$. 
\end{conjecture}


\def\cprime{$'$} \def\cprime{$'$} \def\cprime{$'$}
\providecommand{\bysame}{\leavevmode\hbox to3em{\hrulefill}\thinspace}
\providecommand{\href}[2]{#2}

\end{document}